\documentclass{article}
\usepackage[utf8]{inputenc}
\usepackage{xcolor,graphicx,float}
\usepackage{amsfonts}
\usepackage{amsmath}
\usepackage{bbm}
\usepackage[a4paper, total={6in, 8in}]{geometry}

\begin{document}

\title{Networks with degree-degree correlations is a special case of edge-coloured random graphs}

\author{\small{S\'{a}muel$^1$, G. Balogh,Gergely Palla$^{1,2,3}$, and Ivan Kryven$^{4,5}$\footnote{i.kryven@uu.nl}}\\
\footnotesize{$^1$Department of Biological Physics, E\"{o}tv\"{o}s University, H-1117 Budapest, Hungary}\\
\footnotesize{$^2$MTA-ELTE Statistical and Biological Physics Research Group, }\\
\footnotesize{H-1117 Budapest, Pázmány P. stny. 1/A, Hungary,}\\
\footnotesize{$^3$ Digital Health and Data Utilisation Team, Health Services Management Training Centre},\\
\footnotesize{ Faculty of Health and Public Administration, Semmelweis University, Budapest},\\
\footnotesize{$^4$Mathematical Institute, Utrecht University, Utrecht, 3584 CD, The Netherlands},\\
\footnotesize{$^5$Centre for Complex Systems Studies, 3584 CE Utrecht, the Netherlands},
}
\date{}

\maketitle

\begin{abstract}
{
In complex networks the degrees of adjacent nodes may often appear dependent -- which presents a modelling challenge. 
We present a working framework for studying networks with an arbitrary joint distribution for the degrees of adjacent nodes by showing that such networks are a special case of edge-coloured random graphs. We use this mapping to study bond percolation in networks with assortative mixing and show that, unlike in networks with independent degrees, the sizes of connected components may feature unexpected sensitivity to perturbations in the degree distribution. The results also indicate that degree-degree dependencies may feature a vanishing percolation threshold even when the second moment of the degree distribution is finite. These results may be used to design artificial networks that efficiently withstand link failures and indicate possibility of super spreading in networks without clearly distinct hubs.}
{degree correlated networks, coloured random graph, percolation
}
\\

2000 Math Subject Classification: 05C80, 90B15, 82B43
\end{abstract}

\section*{Introduction}

The random graph with an arbitrary degree distribution is a widely discussed network model in which the links between the nodes obey the maximum entropy principle and the degree distribution is fixed as the only input parameter \cite{newman2003structure}. Such a model is often used as a null-model or `baseline' allowing to detect deviations from the maximum entropy principle in empirical networks and thus presence of useful information \cite{newman2003structure}.

One trivial consequence of the maximum entropy principle is that the degrees of the nodes in such models are independent random variables, whereas in many empirical networks this is not believed to be the case, and significant mutual dependences between the degrees of adjacent nodes are frequently reported \cite{Litvak2013,estrada2011combinatorial}.  A paradigmatic example of positive degree-degree correlation is a network of social contacts where hubs, {\it i.e.} nodes with many neighbours connect to other hubs more frequently than a random choice would dictate \cite{newman2002assortative}.  Conversely, an optimal packing of hard spheres in a finite volume, for example,  is expected to feature negative degree-degree correlation in the corresponding  network of the {\it ad hoc} contacts  \cite{papadopoulos2018network,torres2018modeling}.
Networks with degree-degree dependence are discussed in the context  World Wide Web \cite{olvera2019pagerank}, co-authorship networks  \cite{mondragon2020estimating}, neural networks  \cite{teller2014emergence}, and dynamical processes on networks  \cite{bogua2003epidemic}.

Other features that are not be attributed to the classical infinite configuration model but are commonly observed in empirical networks include clustering~\cite{watts1998collective}, cliques~\cite{palla2005uncovering}, small cycles, multiple edges~\cite{kivela2014multilayer}, hidden embedding space~\cite{serrano2008self}, modular structure~\cite{fortunato2010community} and edges of multiple types or colours. The latter case has been recently addressed in the random graph model with coloured edges \cite{kryven2019bond}. In this model the degree of a node is not a scalar quantity but a vector counting number of edges of each colour and therefore the model is defined by a multivariate degree distribution. It turns out that having freedom to choose such a multivariate distribution may allow one to impose  structural constrains on the network by solely manipulating the degree distribution. For example, specifically chosen distribution may constrain a network to be modular \cite{kryven2019bond} or to have a predefined volume growth trend \cite{Schamboeck2020a}. Even in the case of a simple directed random graph, which can be regarded as random graph with two types of half edges, manipulating the bivariate degree distribution counting in- and out-edges allows to manipulate the location of the percolation transition  \cite{Kryven2016,kryven2017finite,schamboeck2019dynamic}. Interdependent percolation in multiplex networks was shown to correspond to edge percolation in branching cell complexes \cite{Ginestra2019}. 

 In this paper, we explore the connection between random graphs with arbitrary degree-degree distributions and the random graphs with coloured edges and show that the former is a special case of the latter. In other words: if multiple types of edges are permitted in the maximum entropy random graph, the edges can absorb structural information about the degree-degree correlations and thus allow us defining a random graph with a given joint degree-degree distribution. We exploit this connection to study the sizes of connected components and the location of the percolation transitions that occurs during random removal of edges in degree-degree correlated networks. Such a mapping to edge-coloured graphs employs two-dimensional colour labels. This can be contrasted with multiplex networks where one-dimensional  labels on edges are used to represent a stack of layers. 

The rest of the paper is organised as follows: first, we establish the mapping between random networks with a joint  degree-degree distribution and networks that have edges of different types, which we refer to as the edge colour. Such coloured networks are defined by specifying multivariate degree distributions.  We then establish the analytical expression for the critical  percolation threshold, size of the giant component, and typical sizes of sub-extensive components and illustrate several unexpected qualitative phenomena that emerge during bond percolation as a consequence of degree dependence. Namely, we show that, if strong degree dependence is present, perturbed degree distribution may result in peculiar behaviour of connected components during bond percolation, which is reminiscent to degenerate percolation transitions. Second, we demonstrate that  the percolation threshold may be vanishing in even if the degree distribution has a finite second moment. The later indicates a manner to construct networks without clearly distinct hubs that are nevertheless super robust during random failure of links.

\section{The mapping}

Consider a network model in which at the ends of a uniformly at random chosen edge one finds nodes of degree $j,k>0$ with probability $e_{j,k}=e_{k,j},$ with $\sum_{j,k=1}^N e_{k,j}=1$, where $N$ is the maximum degree. Therefore $e_{j,k}$ is a bivariate probability mass function providing the only input information to the model. In all other respects the network is regarded as random. We will also use two quantities that are directly related to $e_{j,k}$:
  \begin{equation}\label{eq:pk}
  p_k=\frac{\sum_{j=1}^N e_{j,k}}{k \sum_{j,k=1}^N k^{-1} e_{j,k}}
\end{equation}
  is the degree distribution, that is the probability that a randomly chosen node has degree $k$, 
  and $$p_{j|k}=\frac{e_{j,k}}{\sum_{j=1}^N e_{j,k}},$$
   is the probability that a node at the end of a randomly chosen edge has degree $j$ given that the node on the other side has degree $k$.

 We will now present a mapping between networks with arbitrary $e_{j,k}$ and a coloured random graph.
In the coloured random graph, each edge is assigned one of $n$ colours, so that a randomly chosen node bears $c_1$ edges of colour one, $c_2$ edges of colour two, and so on. 
Thus the coloured degree of each node can be described by a vector of colour counts denoted by ${\bf c}=\left(c_1,c_2,...,c_n\right)$, and the degree distribution is the  probability  $p(\bf{c})$ that a uniformly at random chosen node has configuration $\bf{c}$. 
The actual degree $k$ of a node with configuration $\bf c$ is given by the sum of all colour counts:
$$k=|{\bf c}|:= \sum_{i=1}^n c_i.$$
Here again, $p(\bf c)$ provides the only information about the model.

Let $N$ be the maximum degree in the network. We consider an edge colouring in which an edge colour $i \in \{1,\dots,N^2\}$ encodes the degrees $j$ and $k$ of the incident nodes to this edge, as given by  the following mapping:
\begin{equation}\label{eq:jkmap}
(j,k)\to i,
\end{equation}
where
$$
\begin{aligned}
&j = (i-1) \text{ div } N +1,\\
&k = (i-1) \text{ mod } N+1.
\end{aligned}
$$
Therefore, we say that a node has $c_i$ edges of colour $i$ if, in accordance with \eqref{eq:jkmap}, it has degree $k$ and it is connected to exactly $c_i$ nodes of degree $j$. As an alternative notation, we write $c_{j,k}:=c_i$, where the mapping between $(j,k)$ and $i$ is understood to be given by \eqref{eq:jkmap}.
 By using the above notations we write the probability that a randomly selected node has adjacent edges of configuration ${\bf c}$ (i.e. $c_1$ edges of colour $1$, $c_2$ edges of color $2$, \emph{etc.}) as a multinomial:
\begin{equation}\label{eq:config_density}
\begin{aligned}
   p({\bf c})=&\mathbb{P}[c_1,\dots,c_{N^2}]=
    \begin{cases}
0, & \exists i\notin \Omega(|{\bf c }|), \text{ and } c_i>0, \\
 p_k\frac{|\bf c|!}{\prod\limits^{N^2}_{i=1}{c_i!}}\ \prod\limits^{N^2}_{i=1}{p^{c_i}_{j|k}}, & \forall i\notin \Omega(|{\bf c }|), \text{ and } c_i=0,
\end{cases}
\end{aligned}
\end{equation}
where  
\begin{equation}
    \Omega(k): = \{i: i\geq (k-1) N \text{ and  }  i< k N \} 
\end{equation}
is the set of colours that may reside on a node of degree $k$.

 The expectations of $c_i$ are given by:
\begin{equation}\label{eq:M1}
\begin{aligned}
\mathbb{E}[c_{j,k}]=&\mathbb{E}[c_i]=\sum_{{\bf c}} c_i p({\bf c}) 
= \sum_{ |{\bf c} | = k}c_{i}  p_k \ \frac{k!}{\prod\limits^{N^2}_{i=1}{c_i!}} 
\ \prod\limits^{N^2}_{i=1} p^{c_i}_{j|k} =k p_k p_{j|k},
\end{aligned}
\end{equation}
for  $j,k=1,\dots,N,$
and to compute the second moments,  $\mathbb{E}[c_{i_1} c_{i_2}]=\mathbb{E}[c_{j_1,k_1} c_{j_2,k_2}]$, we  distinguish three cases:
\begin{enumerate}
\item If  $k_1\neq k_2,$ then
$\mathbb{E}[c_{j_1,k_1}  c_{j_2,k_2}]=0$.
\item  If $j_1\neq j_2$ and $k_1=k_2=k,$ then
\begin{equation*}
\begin{aligned}
\mathbb{E}[c_{j_1,k_1}  c_{j_2,k_2}] = k p_k\left(kp_{j_1|k}p_{j_2|k}-p_{j_1|k}p_{j_2|k}\right).
 \end{aligned}
\end{equation*}\item  If $i=i_1= i_2$ and $k_1=k_2 =k$ then
\begin{equation*}
\mathbb{E}\left[c^2_{j,k}\right]=k p_k \left(k p^2_{j|k}+p_{j|k}(1-p_{j|k})\right).
\end{equation*}
\end{enumerate}
Combining the above-stated cases together gives:
\begin{equation}\label{eq:M2}
\begin{aligned}
&\mathbb{E}[c_{i_1}  c_{i_2}] =\mathbb{E}[c_{j_1,k_1}  c_{j_2,k_2}] = 
\delta_{k_1,k_2} k_1p_{k_1} p_{j_1|k_1}   \left(   ( k_1 - 1 )  p_{j_2|k_1} +  \delta_{j_1,j_2}\right), \; j_1,j_2,k_1,k_2=1,\dots,N
\end{aligned}
\end{equation}

\section{Size of the giant component}
Let  ${\bf P}$ be a permutation matrix with all elements zero except $P_{i_1,i_2}=1$ when colour $i_1$ is identified with colour $i_2$, that is when
$i_1=(j,k)$ and $i_2=(k,j)$. By using the multi-index notation this is written as: $$P_{(j_1,k_1),(j_2,k_2)}=\delta_{j_1,k_2}\delta_{j_2,k_1}.$$
The size of a giant  component in a coloured directed network  \cite{kryven2019bond} is given by: 
\begin{equation}\label{eq:g1}
s= 1- \mathbb{E}[ \bf P x^{c} ],
\end{equation}
where ${\bf x}=\{x_1,\dots,x_{N^2} \}^{\top},\; x_i\in(0,1],$ is the solution of the system 
\begin{equation}\label{eq:g2}
\bf x =  P  F(x),
\end{equation}
with
$$
F({\bf x})_i=\frac{\mathbb{E}[c_i {\bf  x}^{{\bf c - e}_i}]}{\mathbb{E}[c_i]},\; i=1,\dots,N^2
$$
and ${\bf e}_i$ being standard basis vectors of size $N^2.$
By expanding the expectation values used in the above-introduced equations
 for our particular choice of the coloured degree distribution we find that
    $$\mathbb{E}[ {\bf x^{c}} ]=\sum^N_{k=1} p_k\left(\sum^{N}_{j=1}p_{j|k} x_{j,k}\right)^k$$
    and
$$
\begin{aligned}
&   \mathbb{E}[c_{j,k} {\bf s}^{{\bf c - e}_{j,k}}] =\frac{1}{\mathbb{E}[c_{j,k}]} \frac{\partial}{\partial x_{j,k}}\mathbb{E}[ {\bf x^{c}} ]=
\frac{1}{\mathbb{E}[c_{j,k}]}\frac{\partial}{\partial x_{j,k}}\sum^N_{k_1=1} p_{k_1}  \left(\sum^{N}_{j_1=1}p_{j_1|k_1} x_{j_1,k_1}\right)^{k_1} =\\
&\frac{1}{\mathbb{E}[c_{j,k}]} k p_{j|k}p_{k}  \left(\sum^{N}_{j_1=1}p_{j_1|k} x_{j_1,k}\right)^{k-1}=
 \left(\sum^{N}_{j_1=1}p_{j_1|k} x_{j_1,k}\right)^{k-1},
\end{aligned}
$$
which allows us to replace \eqref{eq:g1}-\eqref{eq:g2} with:
\begin{equation}\label{eq:g3}
\begin{aligned}
s =&1- \sum^N_{k=0} p_k\left(\sum^{N}_{j=1}p_{j|k} x_{j,k}\right)^k,\\
x_{k,j}=&\left(\sum^{N}_{l=1}p_{l|k} x_{l,k}\right)^{k-1},\; j,k=1,\dots,N.
\end{aligned}
\end{equation}
Since the right hand side of the latter equation does not depend on $j$, we conclude that: $$x_{k,1}=x_{k,2}=\dots=x_{k,N}.$$
Let $y_{j,k}^{k-1}:=x_{j,k}$,  we may then rewrite \eqref{eq:g3} as
\begin{equation}\label{eq:g}
s =1- \sum^N_{k=0} p_k y_k^k, \;y_{k}=\sum^{N}_{j=1}p_{j|k} y_{j}^{j-1},\; y_{k}\in(0,1].
\end{equation}
Therefore, we have expressed the size of the giant component $s$ in terms the solution of a system with $N$ non-linear equations. Equation \eqref{eq:g} was first presented in  \cite{newman2002assortative} without derivation.

The expression for the expected size of a sub-extensive connected component in coloured random graphs is given by \cite{kryven2019bond}:
\begin{equation}\label{eq:H}
w=\frac{ {\bf x  D(I-H(x) P)}^{-1} \bf{x}}{1-s}+1,
\end{equation}
where $${\bf D}=\text{diag}\{\mathbb E[c_1], \mathbb E[c_2],\dots, \mathbb E[c_{N^2}]\}$$ and $\bf H$ has elements:
$$
\begin{aligned}
&H_{i_1,i_2}({\bf x})=H_{(j_1,k_1),(j_2,k_2)}({\bf x})=
\frac{\mathbb{E}[(c_{i_1}c_{i_2}-\delta_{i_1,i_2}c_{i_1}) {\bf  x}^{{\bf c }- {\bf e}_{i_1} -{\bf e}_{i_2}}]}{\mathbb{E}[c_i]}=\\
&\frac{1}{\mathbb{E}[c_{j_1,k_1}]} \frac{\partial^2}{\partial x_{j_1,k_1}\partial x_{j_2,k_2}}\mathbb{E}[ {\bf x^{c}} ]=
\frac{\partial}{\partial x_{j_2,k_2}}\left(\sum^{N}_{l=1}p_{l|k_1} x_{l,k_1}\right)^{k_1-1}=\\
&(k_1-1)\left(\sum^{N}_{l=1}p_{l|k_1} x_{l,k_1}\right)^{k_1-2}\sum^{N}_{l=1}p_{l|k_1} \delta_{l,j_2}\delta_{k_1,k_2}=
\delta_{k_1,k_2}(k_1-1) p_{j_2|k_1} \left(\sum^{N}_{l=1}p_{l|k_1} x_{l,k_1}\right)^{k_1-2}.\\
\end{aligned}
$$

\section{Bond percolation with degree-degree dependence}
From the theory for edge-coloured random graphs \cite{kryven2019bond}, we know that such networks percolate when 
$$\det({{\bf P}{\bf M}-{\bf I}})=0$$
 where
\begin{equation}\label{eq:M}
M_{i_1,i_2}  = \frac{ \mathbb E[ c_{i_1} c_{i_2} ]}{\mathbb E[ c_{i_2}]}-\delta_{i_1,i_2},\ \; i_1,i_2=1,\dots,N^2,
\end{equation}
and, after plugging the moments expressions \eqref{eq:M1}-\eqref{eq:M2} into \eqref{eq:M},
we obtain
\begin{equation}\label{eq:Mdef}
M_{i_1,i_2} = M_{(j_1,k_1),(j_2,k_2)}= 
 \delta_{k_1,k_2} p_{j_1|k_1} (k_1-1). 
\end{equation}
The elements of the product are given by:
\begin{equation}
({\bf P}{\bf M})_{(j_1,k_1),(j_2,k_2)}=\delta_{j_1,k_2}(k_2-1)p_{k_1|k_2}.
\end{equation}
Even though ${\bf P}{\bf M}$ is a square matrix of size $N^2$,
one can see from the definition \eqref{eq:Mdef}  that this matrix has at most $N$ unique columns and therefore the spectrum of ${\bf P} {\bf M}- {\bf I}$ contains eigenvalue $\lambda =1$ with multiplicity of at least $N^2-N$. 
Since the determinant can be written as the product of all eigenvalues, there exists a smaller matrix of size at most $N$ containing  all the eigenvalues of ${\bf P} {\bf M}- {\bf I}$ apparat of $\lambda=1$ and therefore having the same determinant. Let ${\bf I}$ be $N\times N$ identity matrix, ${\bf e}= (1,1,\dots,1)$ vector of length $N$ and ${\bf S}=N^{-\frac{1}{2}}{\bf I}\otimes {\bf e},$ then 
$$
{\bf C}=  {\bf S}  {\bf P} {\bf M}  {\bf S}^\top,
$$
with
\begin{equation}\label{eq:treshold}
 \det( {{\bf C}- {\bf I}} )= \det( {{\bf P} {\bf M}- {\bf I}} ),
\end{equation}
where $\bf C$  has elements
$$
C_{j,k} = (k-1)p_{j|k},
$$
and size $N$ which is computationally more favourable than size of ${\bf P} {\bf M}$ when detecting phase transitions.

Let us now consider a dynamic network in which  ${\bf C}(\pi)$ is continuously dependent on some parameter $\pi$ and ${\bf C}(0)$ corresponds to no edges, that is $\lim\limits_{\pi \to 0}{\bf C}(\pi)={\bf 0}$. 
In that case $\det( {\bf C}(0) - {\bf I})=\pm 1$ and for some small $\varepsilon$ the determinant does not change its sign when $\pi\in(0,\varepsilon)$, that is $(-1)^{N-1}\det( {\bf C}(\pi) - {\bf I})<0$ for all $\pi<\varepsilon$.
Then the percolation threshold is expressed as
\begin{equation}\label{eq:pc}
\pi_c=\inf\{\pi: (-1)^{N-1}\det({{\bf C}}(\pi_c)-{{\bf I}})>0\},
\end{equation}
that is the smallest $\pi$ for which $(-1)^{N-1}({\bf C}-{{\bf I}})$ has positive determinant. This provides a framework for studying resilience of networks with (dis-)assortative mixing that evolve according to a wide class of dynamic processes. One example of such dynamic process is  bond percolation. 

If $\bf C(\pi)$ is a non-linear function in $\pi$, it is not possible to express the percolation threshold $\pi_c$ as the largest eigenvalue of ${\bf C}$, and one must apply the criticality condition \eqref{eq:pc} to infer the value of the critical percolation threshold. As such, our equation \eqref{eq:pc}, although it is harder to apply, improves previous results on percolation in degree correlated networks \cite{goltsev2008percolation,boguna2002epidemic}, which approximate $ {\bf C}(\pi)$ to be a linear function in $\pi$. Beyond percolation, non-linearity of $ {\bf C}(\pi)$ is typical for dynamic networks that evolve due to some external forcing, for example, as in  \cite{Schamboeck2020b}.

\subsection{Evolution of the joint degree-degree distribution under bond percolation}
Removing edges uniformly at random during bond percolation affects the degree-degree distribution in the following manner,
\begin{equation}
\begin{aligned}\label{eq:ep}
e_{j,k}(\pi)=&\sum_{j_1=j,k_1=k}^{N^2}e_{j_1,k_1}\binom{j_1-1}{j-1}\pi^{j-1}(1-\pi)^{j_1-j}
  \binom{k_1-1}{j-1}\pi^{k-1}(1-\pi)^{k_1-k}.
\end{aligned}
\end{equation}
The degree distribution $p_k(\pi)$ can be readily expressed from the degree-degree distribution \eqref{eq:ep} by applying \eqref{eq:pk} and additionally taking care of the isolated nodes:
$$p_0= \sum_{k>0} p_k \pi^{k-1}.$$

The parameter dependent degree-degree distribution \eqref{eq:ep} gives rise to ${\bf C}(\pi)$ which is generally a non-linear function of $\pi$, and therefore, one must apply equation \eqref{eq:pc} to detect the percolation threshold.  One peculiar property of $e_{i,j}(\pi)$ is that the joint degrees are becoming less dependent with edge removal. This could be seen by studying the correlation coefficient.

\subsection{Decay of the Pearson correlation coefficient during percolation}
The Pearson correlation coefficient $r$ for adjacent degrees depends on $\pi$:
$$
r(\pi)=\frac{r_1}{1-a(1-\frac{1}{\pi})},
$$
where 
$$r_1=\frac{\mathbb{E}_e[jk]-\mathbb{E}_e^2[k]}{\mathbb{E}_e[k^2]-\mathbb{E}_e^2[k]}$$
is the correlation coefficient at $\pi=1$ and
$$
a=\frac{\mathbb{E}_e[k]-1}{\mathbb{E}_e[k^2]-\mathbb{E}_e^2[k]},
$$
characterises how fast the correlation decays.
Indeed, a simple analysis shows that, as long as $r_1>0$, $r(\pi)$ is a strictly increasing function in $\pi$, which means that  uniform removal of edges will always decrease the correlation between adjacent degrees. Moreover,  $r(\pi)$ is convex for $a>1$, concave for $0<a<1$, and a linear function for $a=1$. Note that the value of $a$ is expressed solely in terms of non-mixed moments, and therefore, it is a property of degree distribution and not the copula that characterises dependency between joint degrees. When  $\mathbb{E}[k^2]=\infty$ and   $\mathbb{E}[k]<\infty$,  the correlation coefficient vanishes $r(\pi) \equiv r_1=0 $ for $\pi>0$, however the degree dependency may be still characterised with different measures \cite{Litvak2013,stegehuis2019degree}.

\begin{figure}[h!]
\begin{center}
\includegraphics[width=\textwidth]{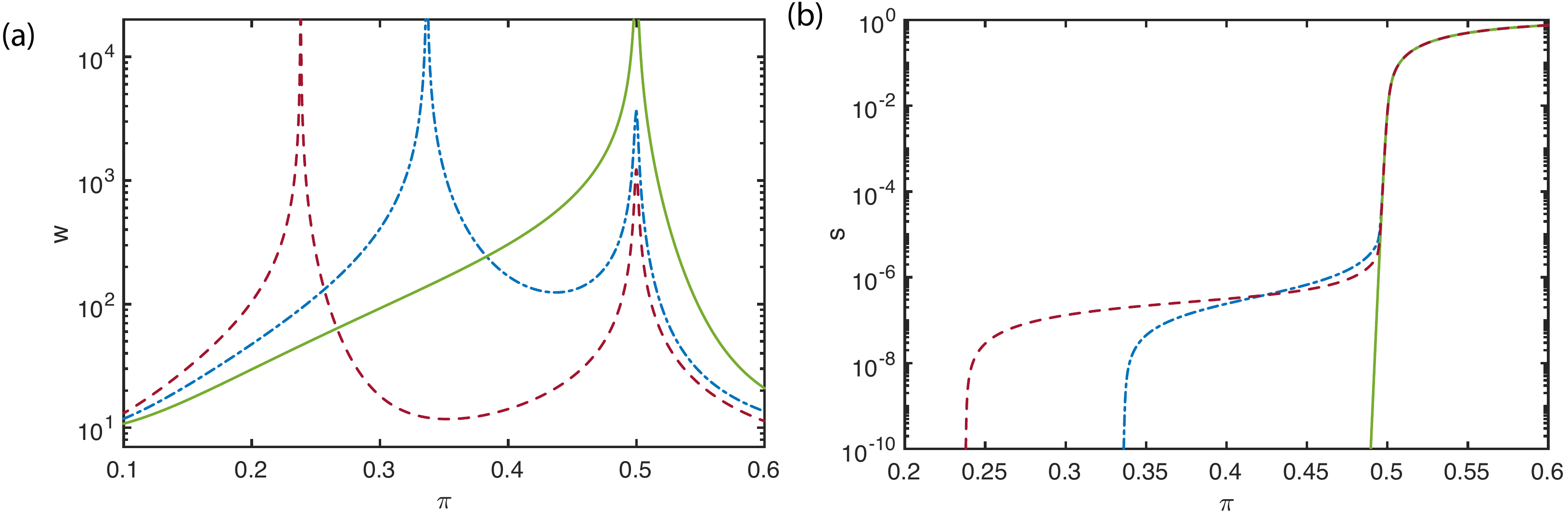}
\caption{
a) Typical size of sub-extensive connected components $w$.
b)  The fraction of nodes in the giant component, $s$.
 The \emph{solid lines} correspond to uncorrected networks with $t=0$,
 \emph{dash-dot lines} to  $t=0.5$, and \emph{dashed lines} to $t=1$.
}
\label{fig:1}
\end{center}
\end{figure}

\begin{figure}[h!]
\begin{center}
\includegraphics[width=0.45\textwidth]{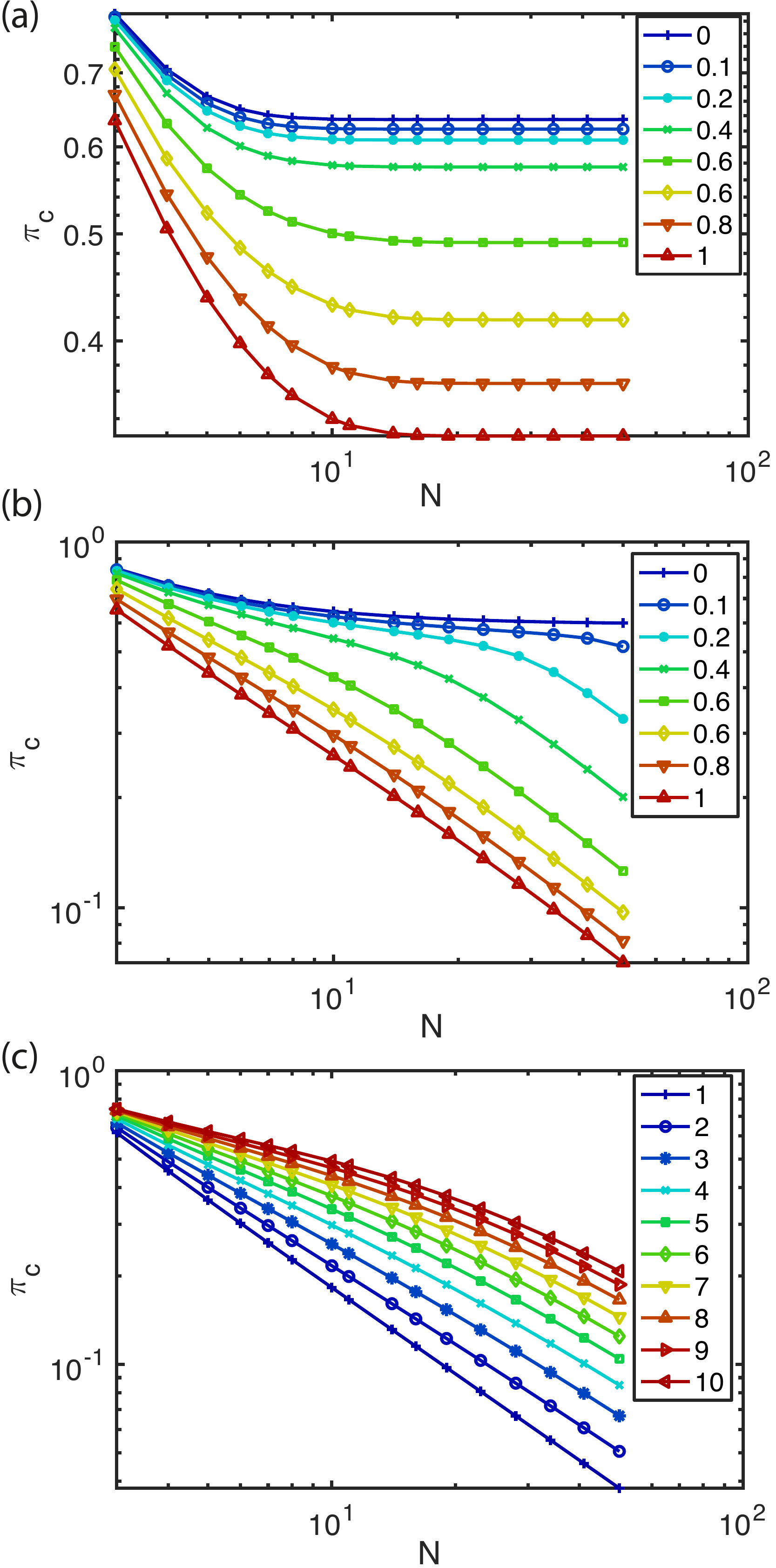}
\caption{
Critical percolation threshold $\pi_c$ plotted versus the maximum degree $N$ for:
1) exponential degree distribution with parameter $t$ as indicated in the legend. 
2) degree distribution with tail exponent $\tau=2.5,$ and parameter $t$ as indicated in the legend.
3) degree distribution with tail exponent $\tau$ as indicated in the legend and  $t=0.9$.
}
\label{fig:2}
\end{center}
\end{figure}

\section{Discussion and conclusions}
 
In this work instead of proposing yet another new model for complex networks we draw the attention to the fact that they can be treated as special cases of coloured random graphs, enabling the collection of many of the existing models under one umbrella. By using this framework, we derive several results for networks defined by their join degree-degree distributions, namely: percolation threshold \eqref{eq:pc}, the size of the giant component \eqref{eq:g}, and typical size of the sub-extensive connected components \eqref{eq:H}. In this section we provide several examples where the link with coloured random graphs reveals unexpected qualitative phenomena in the behaviour of degree-degree correlated networks during random removal of edges.

\subsection{Degenerate percolation transitions}
Consider the degree-degree distribution given by
\begin{equation}\label{eq:ekjkex}
e_{j,k}=(1-t) f(j)f(k)+t \delta_{j,k}f(j),
\end{equation}
where
\begin{equation}\label{eq:bimod}
f(k) = (1-\varepsilon)\delta_{k,3}+\varepsilon\delta_{k,9}.
\end{equation}
with $\varepsilon =10^{-5}$.
Clearly $t>0$ implies dependency between join degrees, with extreme case $t=1$ signifies that the network is composed of multiple disjoint regular graphs.
The parameter $\varepsilon$ is small enough not to induce large changes in the size of the giant component, see Figure~\ref{fig:1}b. 
However, in the correlated case, when $t$ is close to 1, we observe two substantial peaks in the typical sizes of the sub-extensive connected components, while in the uncorrelated case, when $t$ is close to zero, there is only one peak, see Figure~\ref{fig:1}a. Note that the degree distribution is not affected by the value of the coupling parameter $t$ and is bimodal in all three cases  studied in Figure~\ref{fig:1}.
  For all $t,$ there is only one singularity, of the type $O(\frac{1}{\pi-\pi_c})$, while the second peak is bounded. The example illustrates that in networks with degree dependencies, small functional perturbations to the degree distribution may cause large changes in the sizes of connected components.

\subsection{Superrobust networks}
Here we consider a joint degree-degree distribution \eqref{eq:ekjkex} where, as before, parameter $t$ controls degree dependence. We analyse two functional forms of $f_k$: exponential distribution $$f_k =C_1 e^{-k}$$ and distribution with a heavy tail $$f_k=C_2k^{-(\tau+1)},$$ where constants $C_1$ and $C_2$ provide normalisation.
In both cases, we exclude isolated nodes, $f_0:=0,$ and isolated doublets $e_{1,1}:=0$.
Figure \ref{fig:2}a shows that in the exponential case, percolation threshold $\pi_c$  converges to a constant as maximum degree $N$ increases. Here stronger coupling $t$ corresponds to smaller values of the threshold. This is in contrast with Figure \ref{fig:2} b and c where $\pi_c$ is calculated for degree distribution with a heavy tail, showing a steady decrease of $\pi_c$ with increasing maximum degree $N$. This tendency is maintained across different tail exponents and values of coupling constant  $t$.  To date, vanishing percolation threshold has been reported 
only for networks with diverging second moment of the degree distribution, that is $\tau\leq2$, a property which is also frequently  inherited by the spreading processes on such networks \cite{boguna2003absence}.  This example indicates that networks with degree-degree dependencies may feature a vanishing percolation threshold even when the second moment of the degree distribution is finite and therefore the nodes degrees are more homogeneously distributed featuring less pronounced hubs. A direct implication of this observation is that it indicates a way to construct networks with zero percolation threshold that nevertheless do not feature strong degree heterogeneity, and therefore, are robust to random and hub-biased failures \cite{shao2015percolation}. Zero percolation threshold also implies that some spreading processes that have been mapped to percolation, such as susceptible-infected-recovered model with instantaneous transmissions, will feature no epidemic threshold on much wider range of networks that was previously believed \cite{boguna2003absence}.

\subsection*{Acknowledgements}
SB acknowledges kind hospitably of the  Centre for Complex Systems Studies at Utrecht and the support from the Swaantje Mondt PhD Travel Fund.
 GP was partially supported by the Hungarian National Research, Development and Innovation Office (grants no. K128780, NVKP\_16-1-2016-0004) and the Research Excellence Programme of the Ministry for Innovation and Technology in Hungary, within the framework of the Digital Biomarker thematic programme of the Semmelweis University.

\end{document}